\documentclass[12pt]{elsarticle}
\makeatletter
\def\ps@pprintTitle{%
 \let\@oddhead\@empty
 \let\@evenhead\@empty
 \def\@oddfoot{}%
 \let\@evenfoot\@oddfoot}
\makeatother
\usepackage{amsfonts}
\usepackage{latexsym}
\usepackage{amsmath}
\usepackage{amssymb}
\usepackage{amscd}
\usepackage{epsfig}
\usepackage{graphicx}
\usepackage{tkz-berge}
\usetikzlibrary{patterns}
\usetikzlibrary{calc}
\usetikzlibrary{intersections}
\usetikzlibrary{decorations.text}
\usetikzlibrary{decorations.pathreplacing}
\usetikzlibrary{shapes.geometric}
\usetikzlibrary{positioning}
\usepackage{tikz}
\addtocounter{MaxMatrixCols}{4}
\newtheorem{definition}{Definition}[section]

\newtheorem{lemma}[definition]{Lemma}
\newtheorem{theorem}[definition]{Theorem}

\newenvironment{proof}{\par {\sc {\bf Proof.}\hskip 5pt}}{\hfill \qed \par}
\newcommand{\undertilde}[1]{\ensuremath{\mathord{\vtop{\ialign{##\crcr
   $\hfil\displaystyle{#1}\hfil$\crcr\noalign{\kern1.5pt\nointerlineskip}
   $\hfil\tilde{}\hfil$\crcr\noalign{\kern1.5pt}}}}}}
\DeclareMathOperator{\maxcr}{\text{maxcr}}

\begin{document}

\begin{frontmatter}

\title{The maximum crossing number of $C_3 \times C_3$}
\author[mh]{Michael Haythorpe}
\author{Alex Newcombe}
\cortext[mh]{Corresponding author: Michael Haythorpe (michael.haythorpe@flinders.edu.au)}

\begin{abstract}
We determine that the maximum crossing number of $C_3 \times C_3$ is 78, which closes the previously best known range of between 68 and 80.  The proof uses several techniques which may be useful in determining the maximum crossing number of other graphs.
\end{abstract}


\end{frontmatter}

\section{Introduction}

It was shown by Piazza et al \cite{piaz94} that the maximum crossing number of $C_3 \times C_3$ lies somewhere between 68 and 80, inclusive.  We resolve this by showing that the maximum crossing number is exactly 78.  The proof uses a mixture of theoretical arguments as well as exhaustive computer searches.  Although the proofs are ad hoc, we believe that the general techniques could be used to determine the maximum crossing number for other, similarly sized graphs.

\section{Preliminaries}
\begin{figure}[h!]
\centering
\begin{tikzpicture}[thick,scale=0.5]
\draw[fill=white!100,inner sep=0.3pt, minimum width=4pt] (-2,6) node[scale=1]{$(a)$};
\node(1)[circle, draw, fill=black!100,inner sep=0pt, minimum width=4pt] at (0,0){};
\node(2)[circle, draw, fill=black!100,inner sep=0pt, minimum width=4pt] at (2,0){};
\node(3)[circle, draw, fill=black!100,inner sep=0pt, minimum width=4pt] at (1,2){};
\node(4)[circle, draw, fill=black!100,inner sep=0pt, minimum width=4pt] at (6,0){};
\node(5)[circle, draw, fill=black!100,inner sep=0pt, minimum width=4pt] at (8,0){};
\node(6)[circle, draw, fill=black!100,inner sep=0pt, minimum width=4pt] at (7,2){};
\node(7)[circle, draw, fill=black!100,inner sep=0pt, minimum width=4pt] at (3,4){};
\node(8)[circle, draw, fill=black!100,inner sep=0pt, minimum width=4pt] at (5,4){};
\node(9)[circle, draw, fill=black!100,inner sep=0pt, minimum width=4pt] at (4,6){};
\draw (1) -- (2);
\draw (2) -- (3);
\draw (3) -- (1);
\draw (4) -- (5);
\draw (5) -- (6);
\draw (6) -- (4);
\draw (7) -- (8);
\draw (8) -- (9);
\draw (9) -- (7);
\draw (1) to[out=-30,in=210] (4);
\draw (4) -- (7);
\draw (7) -- (1);
\draw (2) to[out=-30,in=210] (5);
\draw (5) -- (8);
\draw (8) -- (2);
\draw (3) -- (6);
\draw (6) -- (9);
\draw (9) -- (3);

\draw[fill=white!100,inner sep=0.3pt, minimum width=4pt] (11,6) node[scale=1]{$(b)$};
\node(10)[circle, draw, fill=black!100,inner sep=0pt, minimum width=4pt] at (15,3){};
\node(11)[circle, draw, fill=black!100,inner sep=0pt, minimum width=4pt] at (13,5){};
\node(12)[circle, draw, fill=black!100,inner sep=0pt, minimum width=4pt] at (13,1){};
\node(13)[circle, draw, fill=black!100,inner sep=0pt, minimum width=4pt] at (17,5){};
\node(14)[circle, draw, fill=black!100,inner sep=0pt, minimum width=4pt] at (17,1){};
\draw (10) -- (11);
\draw (10) -- (12);
\draw (10) -- (13);
\draw (10) -- (14);
\draw (11) -- (12);
\draw (13) -- (14);

\end{tikzpicture}
\caption{The graph $C_3\times C_3$ is shown in $(a)$ and a bowtie graph is shown in $(b)$.\label{figbowtie}}
\end{figure}
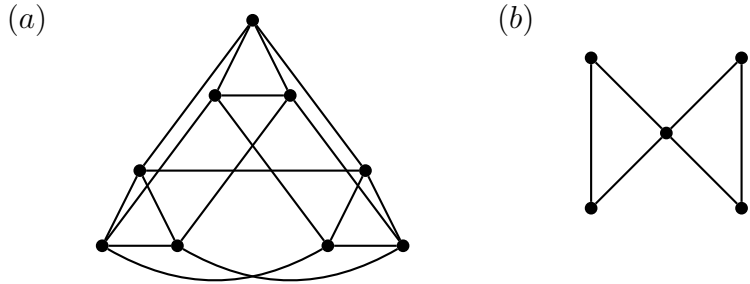

We assume that the reader is familiar with the usual definitions from the theory of crossing numbers and graph drawings.  For more detailed descriptions, we recommend \cite{Schae}.  A graph $G$ has the vertex set $V(G)$ and edge set $E(G)$.  Let $C_n$ denote the cycle graph on $n$ vertices, and let $P_m$ denote the path graph on $m$ edges.  Let $G \times H$ denote the graph Cartesian product between graphs $G$ and $H$.  For example, $C_3 \times C_3$ is displayed in Figure \ref{figbowtie} $(a)$.  A {\em drawing}, $D(G)$, provides a mapping of $E(G)$ and $V(G)$ into the plane.  Vertices are mapped to distinct points and each edge $e=\{u,v\}$ is mapped to a conitnuous arc between the points associated with $u$ and $v$ in such a way that the interior of the arc does not contain any points associated with vertices.  In addition, the interiors of the arcs are only allowed to intersect at singleton points, these are the {\em crossings} and the number of crossings of a drawing is denoted $cr_D(G)$.  We shall refer to the points and arcs given by a drawing as the `vertices' and `edges' of the drawing.  A drawing is a {\em good drawing} if an edge never crosses itself, incident edges do not cross each other and pairs of edges cross at most once.  All drawings considered here are good and when we use the word `drawing', we will mean `good drawing'.  The maximum crossing number, denoted as $\maxcr(G)$, is the maximum number of edge crossings in any good drawing of $G$.   As first described by Conway and Woodall \cite{wood69}, the {\em thrackle number} of $G$, denoted $Th(G)$, is the number 
\begin{equation}
Th(G):=\sum_{\{u,v\} \in E(G)}\frac{1}{2}(|E(G)| - d(v) - d(u) + 1),
\end{equation}
where $d(v)$ is the degree of vertex $v$.  Obviously, for a given graph, the maximum crossing number may not be equal to the thrackle number, but if it is, then $G$ is {\em thrackleable}.  Given a drawing $D$ of a graph $G$, if a pair of non-incident edges do not cross in $D$, then we say that they are a {\em missed pair}.  Thus, for any good drawing $D$, $cr_D(G)$ is equal to the thrackle number minus the number of missed pairs.  Note that the graph $C_4$ is not thrackleable, and similarly, it was shown in \cite{bowtieref} that the {\em bowtie} graph, displayed in Figure \ref{figbowtie} $(b)$, is also not thrackleable.  Thus, any drawing of $C_4$ or a bowtie has at least one missed pair.  The {\em sub-thrackle number}, denoted $STh(G)$, is $Th(G)$ minus the number of subgraphs of $G$ that are isomorphic to $C_4$ plus the number subgraphs isomorphic to the complete graph on 4 vertices.  It follows that, for any graph $G$, $\maxcr(G) \leq STh(G)$.  For a graph $G$, we define a {\em prescription}, denoted $P$, to be a mapping $P: E(G) \rightarrow 2^{E(G)}$, where $2^{E(G)}$ is the power set of $E(G)$.  For each edge $e$, $P(e)$ gives an unordered set of edges, which we call the {\em virtual crossings} of $e$.  If a drawing $D$ of $G$ is such that each edge $e$ crosses exactly those edges in $P(e)$, then $D$ {\em satisfies} the prescription $P$. Clearly, any drawing of $G$ satisfies some prescription $P$, however, it may be the case that there is no drawing that satisfies a given $P$.  We shall be concerned with identifying when it is impossible to satisfy $P$.

After discussing some technical lemmas in Section \ref{secprop}, we shall show the following in Section \ref{secproof}

\begin{theorem}
\label{thm1}
The maximum crossing number of $C_3 \times C_3$ is 78.
\end{theorem}

\section{Properties used in the proof of Theorem \ref{thm1}}
\label{secprop}
Given a drawing of a graph, let $m(A,B)$ be the number, modulo 2, of missed pairs with one of the edges belonging to the edge set $A$ and the other edge belonging to edge set $B$. A property that has been used several times before, e.g. see \cite{wood69}, is as follows

\begin{lemma}\label{lem1}
Let $A$ and $B$ be two vertex disjoint subgraphs of a graph $G$, where $A$ is a cycle of length $n$ and $B$ is a cycle of length $m$.  Then, in any drawing of $G$, $m(E(A),E(B)) = nm \pmod{2}$.
\end{lemma}

For a prescription $P$ of $G$, we shall say that $P$ satisfies {\em property 1} if, for the edges of any two vertex disjoint cycles of $G$, the virtual crossings given by $P$ do not contradict the condition in Lemma \ref{lem1}.

\begin{figure}[h!]
\centering
\begin{tikzpicture}[thick,scale=0.4]
\draw[fill=white!100,inner sep=0.3pt, minimum width=4pt] (-18,4) node[scale=1]{$(a)$};
\node(1)[circle, draw, fill=black!100,inner sep=0pt, minimum width=4pt] at (-18,-2){};
\node(2)[circle, draw, fill=black!100,inner sep=0pt, minimum width=4pt] at (-18,2){};
\node(3)[circle, draw, fill=black!100,inner sep=0pt, minimum width=4pt] at (-16,0){};
\node(4)[circle, draw, fill=black!100,inner sep=0pt, minimum width=4pt] at (-10,-2){};
\node(5)[circle, draw, fill=black!100,inner sep=0pt, minimum width=4pt] at (-10,2){};
\node(6)[circle, draw, fill=black!100,inner sep=0pt, minimum width=4pt] at (-12,0){};
\draw (1) -- (2);
\draw (2) -- (3);
\draw (3) -- (1);
\draw (4) -- (5);
\draw (5) -- (6);
\draw (6) -- (4);
\draw (1) -- (4);
\draw (2) -- (5);
\draw (3) -- (6);

\draw[fill=white!100,inner sep=0.3pt, minimum width=4pt] (-5,4) node[scale=1]{$(b)$};
\foreach \a in {7,8,...,12}{
\node(\a)[circle, draw, fill=black!100,inner sep=0pt, minimum width=4pt] at (\a*360/6: 4cm){};
}
\draw (7) -- (9) -- (11) -- (7);
\draw (8) -- (10) -- (12) -- (8);
\draw (7) -- (10); \draw (8) -- (11); \draw (9) -- (12);

\draw[fill=white!100,inner sep=0.3pt, minimum width=4pt] (8,4) node[scale=1]{$(c)$};
\node(1)[circle, draw, fill=black!100,inner sep=0pt, minimum width=4pt] at (10,-3){};
\node(2)[circle, draw, fill=black!100,inner sep=0pt, minimum width=4pt] at (10,3){};
\node(3)[circle, draw, fill=black!100,inner sep=0pt, minimum width=4pt] at (8.5,0){};
\node(4)[circle, draw, fill=black!100,inner sep=0pt, minimum width=4pt] at (15,-3){};
\node(5)[circle, draw, fill=black!100,inner sep=0pt, minimum width=4pt] at (15,3){};
\node(6)[circle, draw, fill=black!100,inner sep=0pt, minimum width=4pt] at (13.5,0){};
\draw (1) -- (2);
\draw (1) -- (4);
\draw (4) -- (5);
\draw (5) -- (1);
\draw (4) -- (3);
\draw (2) to[out=-60,in=230] (15.55,-3.15) to[out=50,in=-40] (6);
\draw (3) to[out=-45,in=230] (15.95,-3.26) to[out=50,in=-20] (6);
\end{tikzpicture}
\caption{The graph $C_3\times P_1$ is shown in $(a)$. $C_3\times P_1$ is drawn with 15 crossings in $(b)$. $DB(5,3,-1)$ is drawn with 10 crossings in $(c)$.\label{figenv}}
\end{figure}
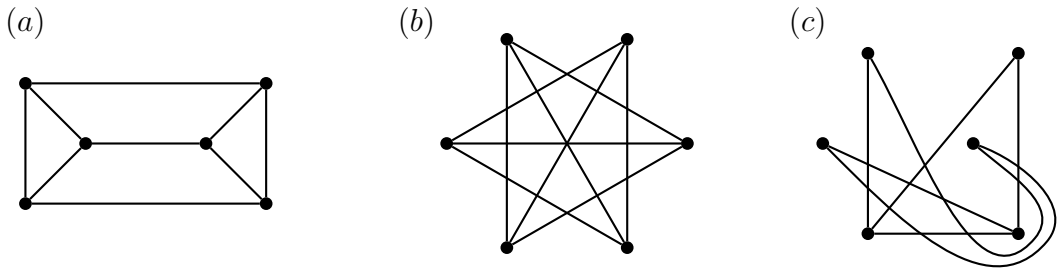

The graph displayed in Figure \ref{figenv} $(a)$ is $C_3 \times P_1$ (sometimes called the envelope graph) and $C_3 \times C_3$ contains 9 subgraphs isomorphic to this graph. We have the following

\begin{lemma}
The maximum crossing number of $C_3 \times P_1$ is 15.
\end{lemma}

\begin{proof}
Figure \ref{figenv} $(b)$ displays $C_3 \times P_1$ drawn with 15 crossings which provides a lower bound and $STh(C_3 \times P_1)=15$ which provides the matching upper bound.
\end{proof}

Next, $C_3\times P_1$ is small enough to use an exhaustive computer search to confirm that there are no drawings of $C_3\times P_1$ with exactly 14 crossings.  Hence, every drawing of $C_3\times P_1$ has either 15 crossings or fewer than 14 crossings.  We note that for larger graphs, exaustive searches such as this quickly become intractable. 

For a prescription $P$ of $G$, we shall say that $P$ satisfies {\em property 2} if, for any subgraph of $G$ which is isomorphic to $C_3 \times P_1$, the virtual crossings given by $P$, restricted to the subgraph, provide a total of either 15 or fewer than 14 virtual crossings.

Next, consider a graph consisting of a cycle of length 6, and one additional edge which forms a triangle.  In the notation of \cite{dumbellref}, this is a $DB(5,3,-1)$ graph.  In \cite{cycletriagref}, Harborth found all thrackleable graphs on six vertices, which does not include $DB(5,3,-1)$ and so there must exist at least one missed pair in any drawing of $DB(5,3,-1)$. Hence, along with the drawing of $DB(5,3,-1)$ in Figure \ref{figenv} $(c)$, which has 10 crossings, we conclude that $\maxcr(DB(5,3,-1))=Th(DB(5,3,-1))-1=10$.

For a prescription $P$ of $G$, we shall say that $P$ satisfies {\em property 3} if, for any subgraph of $G$ which is isomorphic to $DB(5,3,-1)$, the virtual crossings given by $P$, restricted to the subgraph, provide at most 10 virtual crossings.

\section{Proof of Theorem \ref{thm1}}
\label{secproof}
\subsection{Lower bound}
\begin{figure}[h!]
\centering
    \includegraphics[width=0.95\textwidth]{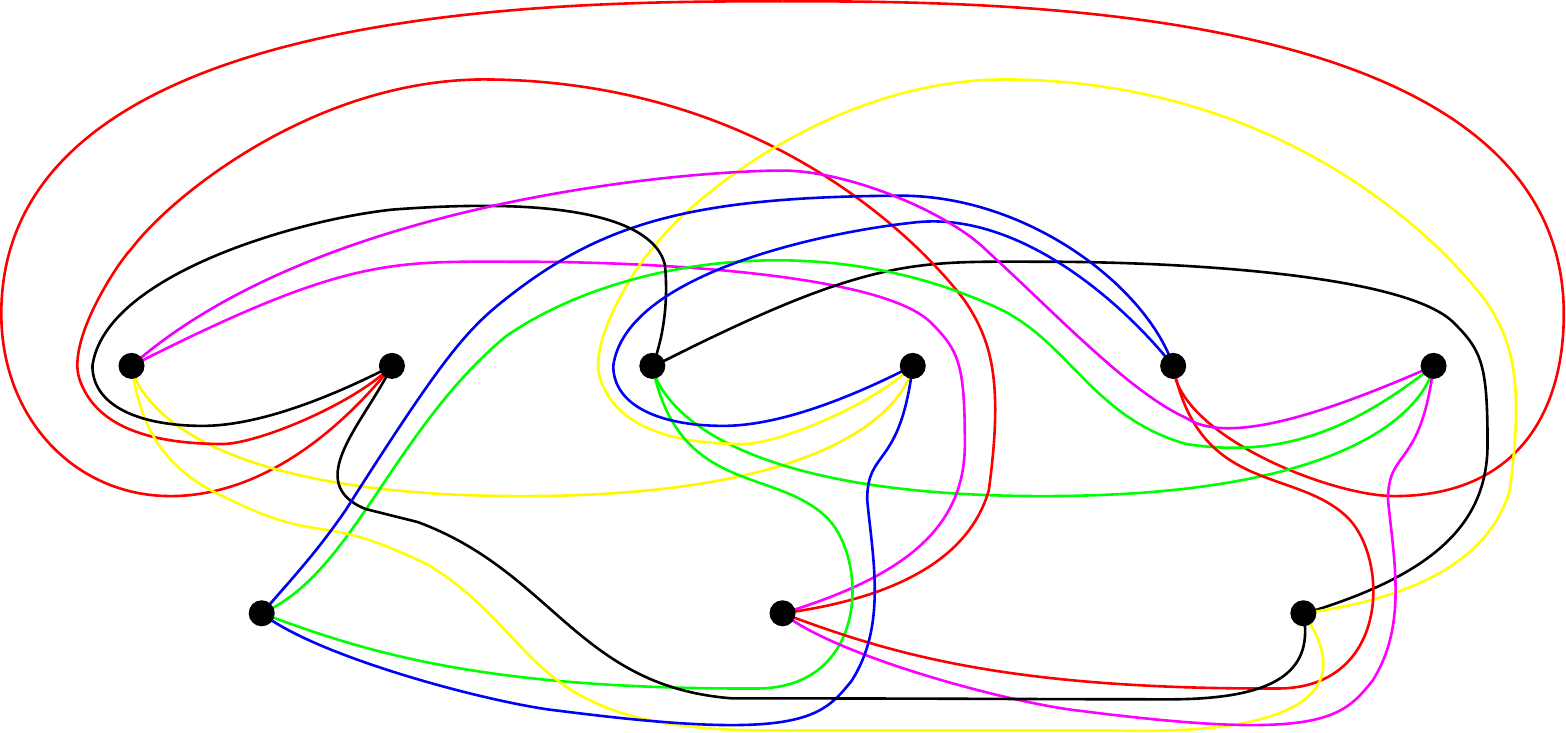}
\caption{$C_3 \times C_3$ drawn with 78 crossings.\label{figc3c3}}
\end{figure}

The lower bound is determined by demonstrating a drawing of $C_3 \times C_3$ with 78 crossings.  The drawing was found using a modified version of the largely successful `planarisation method' for minimising crossings \cite{buch,quickcross}.  The modifications will be discussed in and upcoming publication and they allow us to attempt to maximise crossings by utilising longest paths.  Figure \ref{figc3c3} displays a drawing of $C_3 \times C_3$ with 78 crossings and so $\maxcr(G) \geq 78$. 

\subsection{Upper bound}
First, let $G=C_3 \times C_3$, and let $G-v$ be the graph obtained by deleting vertex $v$ from $G$.  Note that $Th(G-v)=55$.  We will now determine that $\maxcr(G-v)=44$.  This is done in the following way.  There are 5 different 4-cycles along with 4 bowtie subgraphs in $G-v$, each of which must have a missed pair and moreover, no missed pair is double counted.  Hence there must be at least 9 missed pairs in any drawing of $G-v$.

To show that there is not exactly nine missed pairs in $G-v$, we do the following exhaustive search.  We consider all possible combinations of 9 missed pairs, and each combination provides a prescription $P$ of $G-v$.  Note that a prescription $P$ must satisfy properties 1-3 of Section \ref{secprop}, otherwise, there does not exist a drawing of $G$ that satisfies $P$.  The exhaustive search determines that there are no such sets of 9 missed pairs, whose resulting prescription satisfies properties 1-3.  

\begin{figure}[h!]
\centering
    \includegraphics[width=0.5\textwidth]{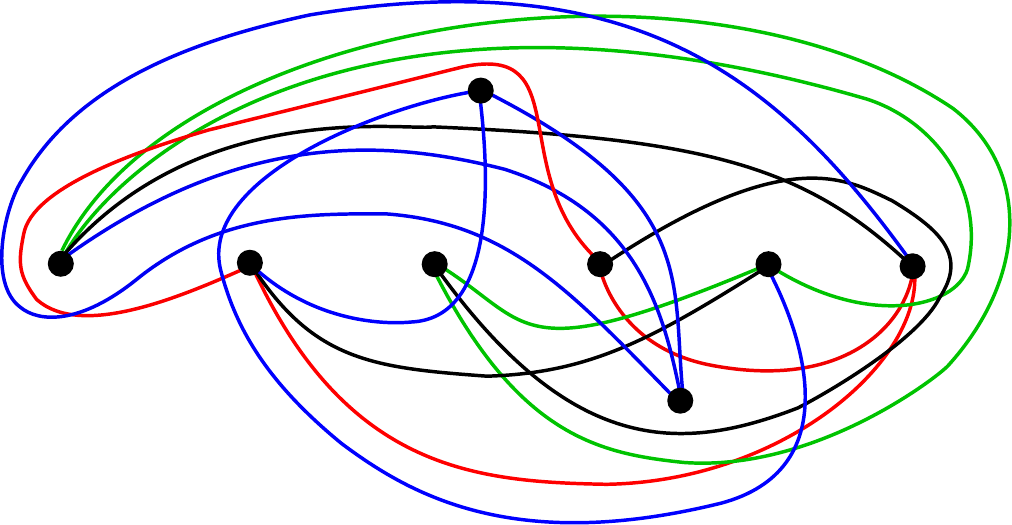}
\caption{$C_3 \times C_3$, minus one vertex, drawn with 44 crossings.\label{figGminusv}}
\end{figure}

To show that there is not exactly 10 missed pairs, we consider all possible combinations of 10 missed pairs in $G-v$ and look for those whose resulting prescription satisfies properties 1-3.  There are 74 such sets, and for each, we have a prescription $P$ for $G-v$.  For each edge $e\in E(G-v)$, $P(e)$ gives a set of edges (the virtual crossings) that $e$ is to cross, however, we do not know the order in which they might be crossed.  We would like to consider all possible orderings and show that it is impossible to produce such a drawing, but this is beyond the realm of tractability.  So instead we do this heuristically in the following way. Consider a subgraph of $G-v$ and the virtual crossings given by $P$, restricted to the edges within the subgraph.  For some subgraphs, it is tractable to consider all possible orderings of the virtual crossings given by $P$.  So, we search for a subgraph for which we are able to verify that it is impossible to draw this subgraph in a way which satisfies the virtual crossings given by $P$, and hence, the supergraph $G-v$ also cannot be drawn to satisfy $P$.  For each of the 74 prescriptions, we were able to find such a subgraph.  At the completion of this process, we have shown that there must be at least 11 missed pairs in any drawing of $G-v$.  Hence, we have that $\maxcr(G-v)\leq 44$ and Figure \ref{figGminusv} displays $G-v$ drawn with exactly 44 crossings which, together, show that $\maxcr(G-v)=44$.

Next, there are 9 possibilites to delete a vertex from $G$ and in any given drawing of $G$, each crossing appears in exactly 5 of the 9 subdrawings corresponding to the subgraph of $G-v$, hence, for any drawing $D$ of $G$ 

\begin{equation}\label{eqD}
cr_D(G)=\frac{1}{5}\sum_{v \in V(G)} cr_D(G-v)\leq \frac{1}{5} (9 \times 44) = 79.2
\end{equation}

and so $\maxcr(G) \leq 79$.  To complete the argument, note that (\ref{eqD}) shows that there is at least one more missed pair than observed by Piazza et al in \cite{piaz94}.  There are two possible cases to consider.  Firstly, if this additional missed pair is formed by edges from two vertex disjoint triangles, then by the arguments in \cite{piaz94}, there is also a second additional missed pair and we would obtain $\maxcr(G) \leq 78$.  

The only other possible case is that the additional missed pair comes from a bowtie subgraph in $G$.  That is, there is a bowtie in $G$ with at least two missed pairs.  Observe that each bowtie of $G$ is contained in 4 of the 9 subgraphs isomorphic to $G-v$.  We now show that if $D$ is a drawing of $G-v$ with two missed pairs on any bowtie, then $cr_D(G-v) < \maxcr(G-v)$.  This again is done by an exhaustive search.  We consider all possible combinations of 11 missed pairs of $G-v$, such that there is a bowtie subgraph with at least 2 of the missed pairs and such that properties 1-3 are still satisfied.  There are 13 such combinations, and for each of these, we find a subgraph that cannot be drawn to satisfy the corresponding prescription.  Hence we have shown that there is no drawing of $G-v$ with $\maxcr(G-v)$ crossings where one of the bowtie subgraphs has at least two missed pairs. That is, in this case,

\begin{equation}
cr_D(G) = \frac{1}{5}\sum_{v \in V(G)} cr_D(G-v)\leq \frac{1}{5} (5 \times 44+4 \times 43) = 78.4.
\end{equation}

We can now conclude that, in either case, $\maxcr(G) \leq 78$, which completes the proof of Theorem \ref{thm1}.

\end{document}